\documentclass[11pt,reqno]{amsart}
\usepackage{amsmath,amssymb,amsfonts,mathrsfs}
\usepackage[margin=1.05in]{geometry}
\usepackage{microtype}
\usepackage{booktabs}
\usepackage[colorlinks=true,linkcolor=blue,citecolor=blue,urlcolor=blue]{hyperref}

\newtheorem{theorem}{Theorem}[section]
\newtheorem{lemma}[theorem]{Lemma}

\newtheorem{corollary}[theorem]{Corollary}
\theoremstyle{definition}
\newtheorem{definition}[theorem]{Definition}
\theoremstyle{remark}
\newtheorem{remark}[theorem]{Remark}
\newtheorem{example}[theorem]{Example}
\numberwithin{equation}{section}

\newcommand{\R}{\mathfrak R}
\newcommand{\I}{\mathfrak I}
\newcommand{\Lh}{\mathscr L}
\newcommand{\Q}{\mathbb Q}
\newcommand{\Z}{\mathbb Z}
\newcommand{\F}{\mathbb F}
\newcommand{\OK}{\mathcal O_K}
\newcommand{\leg}[2]{\left(\frac{#1}{#2}\right)}

\begin{document}

\title[The equation $Cx^2+p_1^{2m_1}\cdots p_k^{2m_k}=y^q$]{On the Diophantine equation $Cx^2+p_1^{2m_1}p_2^{2m_2}\cdots p_k^{2m_k}=y^q$}

\author[U. K. Sangale]{Usha K. Sangale}
\address{School of Mathematical Sciences, Swami Ramanand Teerth Marathwada University, Nanded 431606, India}
\email{ushas073@gmail.com}
\thanks{Corresponding author: Usha K. Sangale (ORCID: \href{https://orcid.org/0000-0003-1033-1594}{0000-0003-1033-1594}).}

\author[P. R. Satpute]{Priyanka R. Satpute}
\address{School of Mathematical Sciences, Swami Ramanand Teerth Marathwada University, Nanded 431606, India}
\email{priyanka.satpute567@gmail.com}


\begin{abstract}
Let $C>3$ be an odd square-free integer, $q$ an odd prime not dividing the class number of $\mathbb Q(\sqrt{-C})$, and $p_1,\dots,p_k$ distinct odd primes not dividing $C$. We describe all solutions of $Cx^2+p_1^{2m_1}\cdots p_k^{2m_k}=y^q$ in positive integers with $\gcd(x,y)=1$ and $y$ odd. Writing $N=p_1^{m_1}\cdots p_k^{m_k}$, every solution comes from a unique factorization $x\sqrt C+N\sqrt{-1}=\pm(u\sqrt C\pm v\sqrt{-1})^q$, or from a half-integral analogue when $q=3$, in which each $p_i\ne q$ divides exactly one of $v$ and $N/v$. Using primitive divisors of Lehmer sequences, we show that for $q\ge5$ some $p_i$ satisfies $p_i\equiv\leg{-C}{p_i}\pmod q$; in particular, there is no solution when $q\ge\max(5,p_1+2,\dots,p_k+2)$. For $q=3$ the equation reduces to four explicit auxiliary equations.
\end{abstract}

\maketitle

%======================================================================
\section{Introduction}\label{sec:intro}
%======================================================================

Equations of the form
\begin{equation}\label{gen}
C_1x^2+C_2=\lambda y^n,\qquad \lambda\in\{1,2,4\},\ n\ge3,
\end{equation}
with fixed positive coefficients $C_1,C_2$, are known as generalized Lebesgue--Ramanujan--Nagell equations; see the survey of Le and Soydan \cite{LS} for background. A recurring theme is the case in which $C_2$ is a product of powers of a few fixed primes. For $\lambda=4$, Chakraborty, Hoque and Srinivas \cite{CHS} solved $Cx^2+p^{2m}=4y^n$, Chakraborty and Hoque \cite{CH} treated $dx^2+p^{2a}q^{2b}=4y^p$, and Hoque, Kotyada and Sangale \cite{HKS} extended this from two prime powers to several, studying
\begin{equation}\label{HKSeq}
cx^2+p_1^{2m_1}p_2^{2m_2}\cdots p_r^{2m_r}=4y^q,\qquad \gcd(x,y)=1,\quad q\nmid h(-c).
\end{equation}
Equations such as $x^2+3^a5^b7^c19^d=4y^n$ were solved completely by Tho \cite{Tho}; see also \cite{BDH,DGS}. For $\lambda=1$ and a single prime, the equation $px^2+q^{2m}=y^p$ was considered by Abu Muriefah \cite{AM08} and Laradji, Mignotte and Tzanakis \cite{LMT}, and the equation
\begin{equation}\label{oneprime}
Cx^2+p^{2m}=y^q
\end{equation}
was classified by Sangale and Satpute \cite{SS} into six explicit families. Patel \cite{Patel} developed a general Lucas--Lehmer approach to \eqref{gen} with $\lambda=1$.

In this paper we generalize the classification of \eqref{oneprime} to an arbitrary number of primes; equivalently, we study the analogue of \eqref{HKSeq} without the factor $4$. The two problems behave differently: in \eqref{HKSeq} the relevant elements of $\OK$ are always half-integral, whereas in our setting they lie in $\Z[\sqrt{-C}]$ except in the single case $q=3$, $C\equiv3\pmod 8$ (Lemma~\ref{replemma}). 

Throughout, we fix an odd square-free integer $C>3$, with $K=\Q(\sqrt{-C})$ and class number $h(-C)$; an odd prime $q$ with $q\nmid h(-C)$; distinct odd primes $p_1,\dots,p_k$ ($k\ge1$) not dividing $C$;
and we consider
\begin{equation}\label{maineq}
Cx^2+p_1^{2m_1}p_2^{2m_2}\cdots p_k^{2m_k}=y^q,\qquad x,y,m_1,\dots,m_k\ge1,\quad \gcd(x,y)=1,\quad y\ \text{odd}.
\end{equation}
We always write
\[
N=p_1^{m_1}p_2^{m_2}\cdots p_k^{m_k},\qquad N=q^{m_0}N_0,\quad q\nmid N_0,
\]
so that $m_0=m_i$ if $q=p_i$ and $m_0=0$ if $q\notin\{p_1,\dots,p_k\}$. The hypothesis that $y$ is odd is automatic unless $C\equiv7\pmod8$ (Lemma~\ref{yodd}).

For integers $u,v$ put
\begin{align}
\R(C,u,v,q)
 &=\sum_{j=0}^{(q-1)/2}\binom{q}{2j}
 u^{q-2j-1}C^{(q-1)/2-j}(-v^2)^j, \label{Rdef}\\
\I(C,u,v,q)
 &=\sum_{j=0}^{(q-1)/2}\binom{q}{2j+1}
 u^{q-2j-1}C^{(q-1)/2-j}(-v^2)^j. \label{Idef}
\end{align}
For $q=3$ these are $\R(C,u,v,3)=Cu^2-3v^2$ and $\I(C,u,v,3)=3Cu^2-v^2$.

\begin{definition}\label{admissible}
A positive divisor $v$ of $N$ is \emph{$q$-admissible} if $m_0\ne1$ and
\[
v=q^{\max(m_0-1,\,0)}\,V,\qquad V\mid N_0,\qquad \gcd\bigl(V,N_0/V\bigr)=1 .
\]
In that case we put $W=N_0/V$, so that $N/v=q^{\min(m_0,1)}\,W$. In words: every $p_i\ne q$ divides exactly one of $v$ and $N/v$, and divides it to the full power $p_i^{m_i}$; and if $q\mid N$, then $q^2\mid N$, $q\mid v$ and $q$ divides $N/v$ exactly once.
\end{definition}

\subsection{Results}
Our first result is a complete parametrization.

\begin{theorem}\label{param}
Suppose that $x,y,m_1,\dots,m_k$ satisfy \eqref{maineq}. Then $x$ is even, $y\equiv1\pmod4$, and there are unique positive integers $u,v$ with $\gcd(Cu,v)=1$ and $v$ a $q$-admissible divisor of $N$ such that exactly one of the following holds.
\begin{enumerate}
\item[\textup{(A)}] $u$ is even and
\[
N=v\,\bigl|\I(C,u,v,q)\bigr|,\qquad x=u\,\bigl|\R(C,u,v,q)\bigr|,\qquad y=Cu^2+v^2 .
\]
\item[\textup{(B)}] $q=3$, $C\equiv3\pmod8$, $u$ and $v$ are odd, and
\[
N=\frac{v\,|3Cu^2-v^2|}{8},\qquad x=\frac{u\,|Cu^2-3v^2|}{8},\qquad y=\frac{Cu^2+v^2}{4}.
\]
\end{enumerate}
In particular, if $q\mid N$ then $q^2\mid N$.

Conversely, let $u,v$ be positive integers with $\gcd(Cu,v)=1$, and suppose that either $u$ is even and $v$ is odd, or $q=3$, $C\equiv3\pmod8$, $u$ and $v$ are odd. Define $N,x,y$ by \textup{(A)} in the first case and by \textup{(B)} in the second, and assume in the second case that $N$ is odd. Then $N,x,y$ are positive integers, $N$ is odd and prime to $C$, $y$ is odd, $\gcd(x,y)=1$, $Cx^2+N^2=y^q$, and $v$ is a $q$-admissible divisor of $N$. The converse does not use the hypothesis $q\nmid h(-C)$.
\end{theorem}

Thus the solutions of \eqref{maineq} correspond bijectively to the pairs $(u,v)$ of Theorem~\ref{param} for which $N$ has the prescribed prime factors. The rigidity of $v$ is the main new phenomenon when $k\ge2$: the representation can distribute the primes $p_i$ between $v$ and $N/v$ only in the ways allowed by Definition~\ref{admissible}. Our second result restricts these distributions.

\begin{theorem}\label{primdiv}
Let $x,y,m_1,\dots,m_k$ satisfy \eqref{maineq}, and let $u,v$ be as in Theorem~\ref{param}.
\begin{enumerate}
\item[\textup{(a)}] Every prime $p_i$ with $p_i\mid N/v$ and $p_i\nmid v$ satisfies $p_i\ne q$ and
\[
p_i\equiv\leg{-C}{p_i}\pmod q ;
\]
in particular $p_i\equiv\pm1\pmod q$.
\item[\textup{(b)}] If $q\ge5$, then \textup{(A)} holds and at least one $p_i$ does not divide $v$; equivalently, $V\ne N_0$ in Definition~\ref{admissible}.
\end{enumerate}
\end{theorem}

Part (a) is elementary and holds for every $q$. Part (b) is where the primitive divisor theorem of Bilu, Hanrot and Voutier \cite{BHV} is used, and only for $q\ge7$; for $q=5$ we use Lucas and Fibonacci numbers. Combining the two parts:

\begin{corollary}\label{sieve}
Let $q\ge5$. If $p_i\not\equiv\leg{-C}{p_i}\pmod q$ for every $i$, then \eqref{maineq} has no solution. In particular \eqref{maineq} has no solution for any prime $q\ge\max(5,\,p_1+2,\dots,p_k+2)$ with $q\nmid h(-C)$. Moreover, in every solution with $q\ge5$, each $p_i\ne q$ with $p_i\not\equiv\leg{-C}{p_i}\pmod q$ satisfies $p_i^{m_i}\mid v$.
\end{corollary}

\begin{example}\label{ex11}
Let $C=11$, so that $h(-11)=1$ and, since $C\equiv3\pmod 8$, $y$ is automatically odd by Lemma~\ref{yodd}. None of $3,5,7$ is congruent to $\pm1$ modulo a prime $q\ge5$, because no such $q$ divides $2,4,6$ or $8$. Hence, by Corollary~\ref{sieve},
\[
11x^2+3^{2a}5^{2b}7^{2c}=y^q,\qquad \gcd(x,y)=1,\ a,b,c\ge1,
\]
has no solution for any prime $q\ge5$.
\end{example}

For $q=3$, Theorem~\ref{param} reduces \eqref{maineq} to four auxiliary equations.

\begin{theorem}\label{qthree}
Let $q=3$, so that $N=3^{m_0}N_0$ with $3\nmid N_0$. If $x,y,m_1,\dots,m_k$ satisfy \eqref{maineq}, then $m_0\ne1$, and there are a positive integer $u$, coprime positive integers $V,W$ with $VW=N_0$ and $\gcd(Cu,V)=1$, and a sign $\varepsilon\in\{\pm1\}$, such that exactly one of the following holds.
\begin{enumerate}
\item[\textup{(A0)}] $m_0=0$, $u$ even, $\varepsilon W\equiv-1\pmod{12}$, and
\[
3Cu^2=V^2+\varepsilon W,\qquad x=u\,|Cu^2-3V^2|,\qquad y=Cu^2+V^2 .
\]
\item[\textup{(A3)}] $m_0\ge2$, $u$ even, $\varepsilon W\equiv1\pmod4$, and
\[
Cu^2=3^{2m_0-3}V^2+\varepsilon W,\qquad x=u\,|Cu^2-3^{2m_0-1}V^2|,\qquad y=Cu^2+3^{2m_0-2}V^2 .
\]
\item[\textup{(B0)}] $m_0=0$, $C\equiv3\pmod8$, $u,V$ odd, $\varepsilon W\equiv1\pmod3$, and
\[
3Cu^2=V^2+8\varepsilon W,\qquad x=\frac{u\,|Cu^2-3V^2|}{8},\qquad y=\frac{Cu^2+V^2}{4} .
\]
\item[\textup{(B3)}] $m_0\ge2$, $C\equiv3\pmod8$, $u,V$ odd, and
\[
Cu^2=3^{2m_0-3}V^2+8\varepsilon W,\qquad x=\frac{u\,|Cu^2-3^{2m_0-1}V^2|}{8},\qquad y=\frac{Cu^2+3^{2m_0-2}V^2}{4} .
\]
\end{enumerate}
Moreover, every prime factor $\ell$ of $W$ satisfies $\ell\equiv\leg{-C}{\ell}\pmod 3$.

Conversely, let $m_0=0$ or $m_0\ge2$, let $u,V,W$ be positive integers with $V$ odd, $\gcd(W,6)=1$ and $\gcd(Cu,V)=1$, and let $\varepsilon\in\{\pm1\}$. If the conditions on $m_0$, $C$, $u$ and $V$ and the displayed equation of one of \textup{(A0)}--\textup{(B3)} hold, then the displayed $x,y$ are positive integers with $Cx^2+(3^{m_0}VW)^2=y^3$, $\gcd(x,y)=1$, $y$ odd and $\gcd(3^{m_0}VW,C)=1$.
\end{theorem}

When one of $V$, $W$ equals $1$, the auxiliary equations become explicit and generalize the six families of \cite{SS}.

\begin{corollary}\label{pure}
Under the hypotheses of Theorem~\ref{param}, consider the following statements about a positive integer $u$.
\begin{enumerate}
\item[\textup{(I)}] $u$ is even and $N=|\I(C,u,1,q)|$, $x=u|\R(C,u,1,q)|$, $y=Cu^2+1$.
\item[\textup{(II)}] $q=3$, $3\nmid N$, $u$ is even and $N^2=3Cu^2+1$, $x=u(8Cu^2+3)$, $y=4Cu^2+1$.
\item[\textup{(III)}] $q=3$, $9\mid N$, $u$ is even and $N^2=27(Cu^2-1)$, $x=u(8Cu^2-9)$, $y=4Cu^2-3$.
\item[\textup{(I$'$)}] $q=3$, $C\equiv3\pmod8$, $u$ is odd and $8N=3Cu^2-1$, $8x=u(Cu^2-3)$, $4y=Cu^2+1$.
\item[\textup{(II$'$)}] $q=3$, $C\equiv3\pmod8$, $3\nmid N$, $u$ is odd and $N^2=3Cu^2-8$, $x=u(Cu^2-3)$, $y=Cu^2-2$.
\item[\textup{(III$'$)}] $q=3$, $C\equiv3\pmod8$, $9\mid N$, $u$ is odd, and for some $\varepsilon\in\{\pm1\}$, $N^2=27(Cu^2-8\varepsilon)$, $x=u(Cu^2-9\varepsilon)$, $y=Cu^2-6\varepsilon$.
\end{enumerate}
A solution of \eqref{maineq} satisfies \textup{(I)} or \textup{(I$'$)} exactly when $v=1$, and satisfies \textup{(II)}, \textup{(III)}, \textup{(II$'$)} or \textup{(III$'$)} exactly when $W=1$. Conversely, if $u\ge1$ and an odd integer $N>1$ satisfy one of these statements, then the displayed $x,y$ are positive integers with $Cx^2+N^2=y^q$, $\gcd(x,y)=1$, $y$ odd and $\gcd(N,C)=1$. If $k=1$, every solution of \eqref{maineq} satisfies exactly one of \textup{(I)}--\textup{(III$'$)}, and if moreover $q\ge5$, only \textup{(I)} occurs and $p_1\ne q$.
\end{corollary}

For $k=1$, Corollary~\ref{pure} is the main theorem of \cite{SS}. For $k\ge2$, the solutions with $v>1$ and $W>1$, which we call \emph{mixed}, are genuinely new: they occur for $q=3$ in all four types of Theorem~\ref{qthree}, and also for $q=5$ and $q=7$, for example
\[
7\cdot94^2+(5\cdot59)^2=53^3,\qquad 5\cdot1990^2+(3\cdot281)^2=29^5,\qquad 5\cdot32006^2+(3\cdot36709)^2=29^7 ,
\]
with $v=5,3,3$, respectively (see Section~\ref{sec:examples}). In the last two examples the Lehmer number attached to the solution has the primitive divisor $281$, respectively $36709$, and Theorem~\ref{primdiv}(a) is visible: $281\equiv1\pmod5$ and $36709\equiv1\pmod 7$.

\subsection*{Organization}
Section~\ref{sec:prelim} contains the representation lemma and the auxiliary results. In Section~\ref{sec:rigid} we prove the rigidity of $v$, and in Section~\ref{sec:param} we prove Theorem~\ref{param}. Section~\ref{sec:primdiv} contains the proof of Theorem~\ref{primdiv} and Corollary~\ref{sieve}, and Section~\ref{sec:three} the proofs of Theorem~\ref{qthree} and Corollary~\ref{pure}. Section~\ref{sec:examples} gives examples and a computational check, and Section~\ref{sec:final} contains concluding remarks.

%======================================================================
\section{Preliminaries}\label{sec:prelim}
%======================================================================

\subsection{The polynomials \texorpdfstring{$\R$}{R} and \texorpdfstring{$\I$}{I}}
The polynomials $\R(C,u,v,q)$ and $\I(C,u,v,q)$ depend only on $u^2$ and $v^2$, and for $\mu\in\{\pm1\}$ the binomial theorem gives
\begin{equation}\label{powerexp}
 (u\sqrt C+\mu vi)^q
 =u\,\R(C,u,v,q)\sqrt C+\mu v\,\I(C,u,v,q)\,i .
\end{equation}
Taking absolute values,
\begin{equation}\label{normid}
C\bigl(u\R(C,u,v,q)\bigr)^2+\bigl(v\I(C,u,v,q)\bigr)^2=(Cu^2+v^2)^q .
\end{equation}
Since $v+u\sqrt{-C}=i\,(u\sqrt C-vi)$, where $\sqrt{-C}=i\sqrt C$, and $i^{q}=(-1)^{(q-1)/2}i$, \eqref{powerexp} with $\mu=-1$ gives
\begin{equation}\label{powerexpK}
(v+u\sqrt{-C})^q=(-1)^{(q-1)/2}\bigl(v\,\I(C,u,v,q)+u\,\R(C,u,v,q)\sqrt{-C}\bigr).
\end{equation}
Reading off single terms of \eqref{Rdef} and \eqref{Idef}, we obtain the congruences
\begin{align}
\I(C,u,v,q)&\equiv(-1)^{(q-1)/2}v^{q-1}\pmod{Cu^2}, & \R(C,u,v,q)&\equiv(-1)^{(q-1)/2}q\,v^{q-1}\pmod{Cu^2},\label{congCu}\\
\I(C,u,v,q)&\equiv(-1)^{(q-1)/2}v^{q-1}\pmod{q}, & \I(C,u,v,q)&\equiv q\,u^{q-1}C^{(q-1)/2}\pmod{v^2}.\label{congqv}
\end{align}
Indeed, in \eqref{Rdef} and \eqref{Idef} every term except the one with $j=(q-1)/2$ is divisible by $Cu^2$; every term of \eqref{Idef} with $j\ge1$ is divisible by $v^2$; and $q\mid\binom{q}{2j+1}$ for $2j+1<q$.

\subsection{The parity of \texorpdfstring{$y$}{y}}

\begin{lemma}\label{yodd}
Let $x,y\ge1$ and an odd integer $N$ satisfy $Cx^2+N^2=y^q$ with $\gcd(x,y)=1$. If $y$ is even, then $x$ is odd and $C\equiv7\pmod8$. In particular $y$ is odd whenever $C\not\equiv7\pmod 8$.
\end{lemma}

\begin{proof}
If $y$ is even then $x$ is odd, so $Cx^2+N^2\equiv C+1\pmod 8$, while $8\mid y^q$ because $q\ge3$. Hence $C\equiv7\pmod 8$.
\end{proof}

\subsection{A representation lemma}

\begin{lemma}\label{replemma}
Let $x,Y,Z$ be positive integers such that
\[
Cx^2+Y^2=Z^q,\qquad Z\ \text{odd},\qquad \gcd(Y,Z)=1 .
\]
Then there is $\gamma\in\OK$ with $Y+x\sqrt{-C}=\gamma^q$, and $\gamma$ is unique. Write $\gamma=(s+t\sqrt{-C})/2$ with $s,t\in\Z$, $s\equiv t\pmod 2$. Then $st\ne0$, and there are positive integers $u,v$ with $\gcd(Cu,v)=1$ and signs $\delta,\mu\in\{\pm1\}$ such that one of the following holds:
\begin{enumerate}
\item[(a)] $s,t$ are even, $(u,v)=(|t|/2,|s|/2)$, $x\sqrt C+Yi=\delta(u\sqrt C+\mu vi)^q$ and $Z=Cu^2+v^2$;
\item[(b)] $s,t$ are odd, $(u,v)=(|t|,|s|)$, $q=3$, $C\equiv3\pmod 8$, and
$x\sqrt C+Yi=\delta\bigl((u\sqrt C+\mu vi)/2\bigr)^3$, $Z=(Cu^2+v^2)/4$.
\end{enumerate}
Conversely, any representation as in \textup{(a)} or \textup{(b)} arises in this way from $\gamma$; in particular the case and the pair $(u,v)$ are determined by $(x,Y)$.
\end{lemma}

\begin{proof}
In $\OK$ we have $(Y+x\sqrt{-C})(Y-x\sqrt{-C})=Z^q$. Let $\mathfrak P$ be a prime ideal dividing both factors. Then $\mathfrak P$ contains $Z^q$, so it lies above an odd prime $\ell\mid Z$; and it contains $2Y$, hence $Y$ because $2\notin\mathfrak P$. Thus $\ell\mid\gcd(Y,Z)=1$, which is absurd. So the two principal ideals are coprime, and since their product is the $q$th power of $(Z)$, we have $(Y+x\sqrt{-C})=\mathfrak a^q$ for some ideal $\mathfrak a$. The class of $\mathfrak a$ has order dividing both $q$ and $h(-C)$, so $\mathfrak a=(\gamma_0)$ is principal. As $C>3$, the units of $\OK$ are $\pm1$; since $q$ is odd we may replace $\gamma_0$ by $\pm\gamma_0$ and obtain $Y+x\sqrt{-C}=\gamma^q$. If also $\gamma'^q=Y+x\sqrt{-C}$, then $\gamma'/\gamma\in K$ is a $q$th root of unity, hence equal to $1$; so $\gamma$ is unique.

Write $\gamma=(s+t\sqrt{-C})/2$ as in the statement, so that $Z=N(\gamma)=(s^2+Ct^2)/4$. If $t=0$ then $\gamma^q\in\Z$ and $x=0$; if $s=0$ then $\gamma^q\in\Z\sqrt{-C}$ and $Y=0$. Hence $st\ne0$.

Suppose that $s,t$ are odd. Then $\gamma\notin\Z[\sqrt{-C}]$, so $C\equiv3\pmod 4$. If $C\equiv7\pmod 8$, then $s^2+Ct^2\equiv0\pmod 8$ and $Z$ is even, a contradiction; hence $C\equiv3\pmod8$. Then $2$ is inert in $K$ and $\OK/2\OK\cong\F_4$. Put $\omega=(1+\sqrt{-C})/2$. Since $\sqrt{-C}=2\omega-1$, the image of $\Z[\sqrt{-C}]$ in $\F_4$ is the prime field $\F_2$. On the other hand $\gamma=(s-t)/2+t\omega$ with $t$ odd, so the image $\bar\gamma$ of $\gamma$ in $\F_4$ is not in $\F_2$, and has order $3$ in $\F_4^\times$. But $\gamma^q=Y+x\sqrt{-C}\in\Z[\sqrt{-C}]$, so $\bar\gamma^q\in\F_2^\times=\{1\}$, whence $3\mid q$ and $q=3$.

Conjugating $Y+x\sqrt{-C}=\gamma^q$ and multiplying by $i$, using $i(Y-x\sqrt{-C})=x\sqrt C+Yi$, $i(s-t\sqrt{-C})=t\sqrt C+si$ and $i=(-1)^{(q-1)/2}i^q$, we obtain
\[
x\sqrt C+Yi=(-1)^{(q-1)/2}\Bigl(\frac{t\sqrt C+si}{2}\Bigr)^q .
\]
Since
\[
t\sqrt C+si=\operatorname{sgn}(t)\bigl(|t|\sqrt C+\operatorname{sgn}(st)|s|\,i\bigr)
\]
and $q$ is odd, this is (a) or (b) with $\delta=(-1)^{(q-1)/2}\operatorname{sgn}(t)$ and $\mu=\operatorname{sgn}(st)$.

To see that $\gcd(Cu,v)=1$, let $\ell$ be a prime dividing $Cu$ and $v$. By \eqref{powerexp}, $Y=\pm v\,\I(C,u,v,q)$ in case (a) and $8Y=\pm v\,\I(C,u,v,3)$ in case (b); so $\ell\mid Y$ in case (a), and also in case (b) because $v$ is odd. Moreover $\ell$ divides $Cu^2+v^2$, which is $Z$ in case (a) and $4Z$ in case (b), so $\ell\mid Z$ (in case (b) since $\ell$ is odd). This contradicts $\gcd(Y,Z)=1$.

Finally, suppose that (a) or (b) holds for some $u,v,\delta,\mu$. Replacing $v$ by $\mu v$ in \eqref{powerexpK} and comparing coefficients shows that $\gamma'=\pm(\mu v+u\sqrt{-C})$ in case (a), and $\gamma'=\pm(\mu v+u\sqrt{-C})/2$ in case (b), satisfies $\gamma'^q=Y+x\sqrt{-C}$ for a suitable sign. Hence $\gamma'=\gamma$, and the case and $(u,v)$ are determined by $\gamma$.
\end{proof}

\begin{remark}\label{yuanremark}
Lemma~\ref{replemma} is a version of a standard ideal-factorization argument; compare Bugeaud and Shorey \cite{BS} and Yuan \cite{Yuan}. Case (b) cannot be omitted: for $C=11$,
\[
11\cdot2^2+3^4=5^3,\qquad
3^2+2\sqrt{-11}=\left(\frac{-3+\sqrt{-11}}2\right)^3,
\]
and $3^2+2\sqrt{-11}$ is not a cube in $\Z[\sqrt{-11}]$. This case is invisible to an argument in the order $\Z[\sqrt{-C}]$, whose class number is $h(-4C)=3h(-C)$ when $C\equiv3\pmod 8$.
\end{remark}

\subsection{Lehmer pairs}
A pair $(\alpha,\beta)$ of algebraic integers is a \emph{Lehmer pair} if $(\alpha+\beta)^2$ and $\alpha\beta$ are non-zero coprime rational integers and $\alpha/\beta$ is not a root of unity. For odd $n$ its $n$th \emph{Lehmer number} is
\[
\Lh_n(\alpha,\beta)=\frac{\alpha^n-\beta^n}{\alpha-\beta},
\]
and for even $n$ one divides by $\alpha^2-\beta^2$ instead. A prime $\ell$ is a \emph{primitive divisor} of $\Lh_n(\alpha,\beta)$ if $\ell\mid\Lh_n(\alpha,\beta)$ but $\ell\nmid(\alpha^2-\beta^2)^2\Lh_1(\alpha,\beta)\cdots\Lh_{n-1}(\alpha,\beta)$. A Lehmer pair is \emph{$n$-defective} if $\Lh_n(\alpha,\beta)$ has no primitive divisor. Two Lehmer pairs $(\alpha_1,\beta_1)$, $(\alpha_2,\beta_2)$ are \emph{equivalent} if $\alpha_1/\alpha_2=\beta_1/\beta_2\in\{\pm1,\pm\sqrt{-1}\}$. The parameters of $(\alpha,\beta)$ are $(a,b)=\bigl((\alpha+\beta)^2,(\alpha-\beta)^2\bigr)$, and passing to an equivalent pair changes $(a,b)$ at most to $(-a,-b)$.

\begin{theorem}[Bilu--Hanrot--Voutier \cite{BHV}; Voutier \cite{Voutier95}]\label{bhv}
If $n>30$, there are no $n$-defective Lehmer pairs. If $q$ is a prime with $7\le q\le 29$ and $(\alpha,\beta)$ is $q$-defective, then $q\in\{7,13\}$ and, up to equivalence, $(a,b)$ is one of
\[
(1,-7),\ (1,-19),\ (3,-5),\ (5,-7),\ (13,-3),\ (14,-22)\quad(q=7);\qquad (1,-7)\quad(q=13).
\]
\end{theorem}

This is \cite[Theorem~1.4]{BHV} together with Table~2 of \cite{BHV}. The recent corrections of Voutier \cite{VoutierCorr} concern only the indices $n\in\{3,4,5,6,8,10,12\}$; since we treat $q=3$ and $q=5$ directly, they do not affect our arguments.

\subsection{Lucas and Fibonacci numbers}
Let $(F_n)$ and $(L_n)$ be the Fibonacci and Lucas sequences and $\varphi=(1+\sqrt5)/2$, so that $\varphi^n=(L_n+F_n\sqrt5)/2$. We use the following classical facts.
\begin{enumerate}
\item[(i)] The integer solutions of $X^2-5Y^2=\pm4$ are $(X,Y)=(\pm L_n,\pm F_n)$, $n\ge0$.
\item[(ii)] The units of $\Z[\sqrt5]$ are $\pm(2+\sqrt5)^n$, $n\in\Z$, and $2+\sqrt5=\varphi^3$; hence every solution of $X^2-5Z^2=\pm1$ in positive integers satisfies $X+Z\sqrt5=\varphi^{3n}$ for some $n\ge1$, so that $Z=F_{3n}/2$.
\item[(iii)] (Cohn \cite{Cohn}) The only Lucas numbers that are perfect squares are $L_1=1$ and $L_3=4$.
\item[(iv)] (Keskin and Karaatli \cite{KK}) If $n\ge1$ and $F_n=5z^2$ with $z$ an integer, then $n=5$.
\end{enumerate}

%======================================================================
\section{Rigidity of the divisor \texorpdfstring{$v$}{v}}\label{sec:rigid}
%======================================================================

The following lemma uses only the formulas in (A) and (B), so it applies both to solutions and in the converse direction.

\begin{lemma}\label{rigid}
Let $u,v$ be positive integers with $\gcd(Cu,v)=1$, and put
\[
J=\I(C,u,v,q)\quad\text{in case (A)},\qquad J=\frac{3Cu^2-v^2}{8}\quad\text{in case (B)},
\]
where in case \textup{(B)} we assume $q=3$, $C\equiv3\pmod8$ and $u,v$ odd. Let $N=v|J|$. Then:
\begin{enumerate}
\item[(i)] every prime dividing both $v$ and $J$ equals $q$;
\item[(ii)] if $q\mid v$, then $q$ divides $J$ exactly once;
\item[(iii)] if $q\nmid v$, then $q\nmid J$.
\end{enumerate}
Consequently, if $J\ne0$ and $N$ is odd, then $v$ is a $q$-admissible divisor of $N$, and $N/v=|J|$.
\end{lemma}

\begin{proof}
Case (A). By \eqref{congqv}, $J\equiv q\,u^{q-1}C^{(q-1)/2}\pmod{v^2}$. If a prime $\ell$ divides $v$ and $J$, then $\ell\mid q\,u^{q-1}C^{(q-1)/2}$, and since $\ell\nmid Cu$ we get $\ell=q$; this is (i). If $q\mid v$, then $q^2\mid v^2$ and $q\nmid Cu$, so the $q$-adic valuation of $J$ equals that of $q\,u^{q-1}C^{(q-1)/2}$, namely $1$; this is (ii). If $q\nmid v$, then $J\equiv\pm v^{q-1}\not\equiv0\pmod q$ by \eqref{congqv}; this is (iii).

Case (B). Here $8J=3Cu^2-v^2$. If a prime $\ell$ divides $v$ and $J$, then $\ell$ is odd and $\ell\mid3Cu^2$, so $\ell=3$ because $\ell\nmid Cu$. If $3\mid v$, then $9\mid v^2$ and $3\nmid Cu$, so $3$ divides $3Cu^2-v^2$ exactly once. If $3\nmid v$, then $3Cu^2-v^2\equiv-v^2\not\equiv0\pmod3$.

Now let $N=v|J|$ be odd. For a prime $p\ne q$ dividing $N$, (i) shows that $p$ divides exactly one of $v$ and $J$, so the full power of $p$ in $N$ lies in $v$ or in $J$. If $q\mid N$, then $q\mid v$ by (iii), and by (ii) $q$ divides $J=N/v$ exactly once; hence $m_0\ge2$ and $q^{m_0-1}\,\|\,v$. If $q\nmid N$, then $q\nmid v$. In all cases $v$ is $q$-admissible.
\end{proof}

%======================================================================
\section{Proof of Theorem~\ref{param}}\label{sec:param}
%======================================================================

\subsection*{Parity}
If $p_i\mid y$ for some $i$, then $p_i\mid Cx^2$, so $p_i\mid x$ as $p_i\nmid C$, contradicting $\gcd(x,y)=1$. Hence
\begin{equation}\label{coprime}
\gcd(N,y)=1 .
\end{equation}
Since $y$ and $C$ are odd, $Cx^2=y^q-N^2$ is even, so $x$ is even; then $4\mid Cx^2$ and $y^q\equiv1\pmod4$, and as $q$ is odd, $y\equiv1\pmod4$.

\subsection*{Parametrization}
By \eqref{coprime}, Lemma~\ref{replemma} applies with $Y=N$ and $Z=y$. Let $u,v,\delta,\mu$ be the resulting parameters. Comparing coefficients by \eqref{powerexp},
\begin{align}
x&=\delta u\,\R(C,u,v,q),& N&=\delta\mu v\,\I(C,u,v,q) &&\text{in case (a)},\label{realimag}\\
8x&=\delta u\,(Cu^2-3v^2),& 8N&=\delta\mu v\,(3Cu^2-v^2) &&\text{in case (b)}.\label{halfrealimag}
\end{align}
In case (a), $v$ divides the odd number $N$, so $v$ is odd; since $y=Cu^2+v^2$ is odd, $u$ is even. Thus case (a) gives (A) and case (b) gives (B). Since $N\ne0$ is odd, Lemma~\ref{rigid} shows that $v$ is a $q$-admissible divisor of $N$; in particular, if $q\mid N$ then $m_0\ge2$. Uniqueness of $(u,v)$, and the fact that exactly one of (A), (B) holds, follow from the last statement of Lemma~\ref{replemma}. Indeed, if $u,v$ satisfy the formulas of (A), then $\R(C,u,v,q)\ne0$ and $\I(C,u,v,q)\ne0$, and with $\delta=\operatorname{sgn}\R(C,u,v,q)$ and $\mu=\delta\operatorname{sgn}\I(C,u,v,q)$, \eqref{powerexp} gives $x\sqrt C+Ni=\delta(u\sqrt C+\mu vi)^q$, a representation as in (a); the same choice of signs turns the formulas of (B) into a representation as in (b).

\subsection*{The converse}
Let $u,v$ be as in the converse statement and put $A=Cu^2$.

\emph{Case (A).} Here $u$ is even and $v$ is odd. By \eqref{normid}, $Cx^2+N^2=y^q$. The number $y=A+v^2$ is odd. By \eqref{congCu}, $\I(C,u,v,q)\equiv\pm v^{q-1}\pmod{A}$; as $A$ is even and $\gcd(A,v)=1$, $\I(C,u,v,q)$ is odd and non-zero, so $N$ is a positive odd integer. Also by \eqref{congCu}, $\R(C,u,v,q)\equiv\pm q\,v^{q-1}\pmod A$; if $\R(C,u,v,q)=0$, then $A\mid q$ because $\gcd(A,v)=1$, which is impossible as $A\ge4C\ge20$. Hence $x>0$. If a prime $\ell$ divides $C$ and $N$, then $\ell\nmid v$, so $\ell\mid\I(C,u,v,q)\equiv\pm v^{q-1}\pmod C$, a contradiction; so $\gcd(N,C)=1$.

Suppose that a prime $\ell$ divides $x$ and $y$. Then $\ell$ is odd and $\ell\mid N^2$, so $\ell\mid N$. If $\ell\mid v$, then $\ell\mid y-v^2=A$, contradicting $\gcd(Cu,v)=1$. Otherwise $\ell\mid\I(C,u,v,q)$. Put $\zeta=v+u\sqrt{-C}\in\OK$. By \eqref{powerexpK}, $\zeta^q=\pm\bigl(v\,\I(C,u,v,q)+u\,\R(C,u,v,q)\sqrt{-C}\bigr)$, and $\ell$ divides both $v\,\I(C,u,v,q)$ and $u\,\R(C,u,v,q)=\pm x$. So $\ell$ divides $\zeta^q$ and $\bar\zeta^q$ in $\OK$, and a prime ideal above $\ell$ divides $\zeta$ and $\bar\zeta$, hence $\zeta+\bar\zeta=2v$. Thus $\ell\mid v$, a contradiction. Hence $\gcd(x,y)=1$.

\emph{Case (B).} Here $q=3$, $C\equiv3\pmod8$ and $u,v$ are odd, so $A\equiv3\pmod 8$ and $v^2\equiv1\pmod8$. Hence $A-3v^2\equiv3-3=0$ and $3A-v^2\equiv9-1\equiv0\pmod8$, while $A+v^2\equiv4\pmod8$. So $x$ and $N$ are integers and $y$ is an odd integer. The identity
\[
A(A-3v^2)^2+v^2(3A-v^2)^2=(A+v^2)^3,
\]
divided by $64$, gives $Cx^2+N^2=y^3$. If $A=3v^2$, then $v=1$ (as $\gcd(A,v)=1$) and $C=3$, which is excluded; if $3A=v^2$, then $3\mid v$ and $9\mid3A$, so $3\mid A$, contradicting $\gcd(A,v)=1$. Hence $x,N>0$. Let $\ell$ be a prime dividing $N$; then $\ell$ is odd. If $\ell\mid v$, then $\ell\nmid A$, so $\ell\nmid C$ and $\ell\nmid 4y=A+v^2$. Otherwise $\ell\mid3A-v^2$; then $\ell\mid C$ would give $\ell\mid v$, and $\ell\mid y$ would give $\ell\mid(3A-v^2)+(A+v^2)=4A$, hence $\ell\mid A$ and $\ell\mid v$; both are impossible. So $\gcd(N,C)=\gcd(N,y)=1$, and therefore $\gcd(x,y)=1$, because a common prime factor of $x$ and $y$ divides $N^2$.

In both cases $N$ is odd and $J\ne0$, so $v$ is $q$-admissible by Lemma~\ref{rigid}. None of these verifications uses $q\nmid h(-C)$. This completes the proof of Theorem~\ref{param}. \qed

%======================================================================
\section{Primitive divisors: proof of Theorem~\ref{primdiv}}\label{sec:primdiv}
%======================================================================

\subsection{Proof of part (a)}
Let $\gamma=(s+t\sqrt{-C})/2$ be the element of Lemma~\ref{replemma} with $\gamma^q=N+x\sqrt{-C}$, so that $|s|=2v$ in case (A) and $|s|=v$ in case (B), and $\gamma\bar\gamma=y$. Let $\ell=p_i$ be a prime with $\ell\mid N$ and $\ell\nmid v$. Then $\ell\ne q$ by Lemma~\ref{rigid}(iii), $\ell\nmid y$ by \eqref{coprime}, and $\ell$ is unramified in $K$ because $\ell$ is odd and $\ell\nmid C$. Let $\mathfrak l$ be a prime ideal of $\OK$ above $\ell$; then $\gamma$ and $\bar\gamma$ are units at $\mathfrak l$. Put
\[
\rho=\frac{\gamma}{\bar\gamma}=\frac{\gamma^2}{y}\in K,\qquad \rho\bar\rho=1 .
\]
Since $\gamma^q+\bar\gamma^q=2N\equiv0\pmod{\mathfrak l}$, we have $\rho^q\equiv-1\pmod{\mathfrak l}$, so the order $d$ of $\rho$ in $(\OK/\mathfrak l)^\times$ divides $2q$ but not $q$ (as $\ell$ is odd). Thus $d\in\{2,2q\}$. If $d=2$, then $\rho\equiv-1$, that is, $\mathfrak l$ divides $\gamma+\bar\gamma=s$, so $\ell\mid s$ and $\ell\mid v$, contrary to assumption. Hence $d=2q$.

If $\leg{-C}{\ell}=1$, then $\OK/\mathfrak l\cong\F_\ell$ and $2q\mid\ell-1$. If $\leg{-C}{\ell}=-1$, then $\OK/\mathfrak l\cong\F_{\ell^2}$ and the Frobenius automorphism acts as complex conjugation, so $\rho^\ell\equiv\bar\rho=\rho^{-1}\pmod{\mathfrak l}$ and $2q\mid\ell+1$. In both cases $\ell\equiv\leg{-C}{\ell}\pmod q$. \qed

\begin{remark}
In case (A), a prime $\ell$ as in part (a) is a primitive divisor of the Lehmer number $\Lh_q(\alpha,\beta)$ attached to the pair $\alpha=u\sqrt C+\mu vi$, $\beta=\bar\alpha$ (see the proof of part (b) below), and part (a) is then a special case of Lehmer's law of apparition \cite{Lehmer}. We have included the short direct proof, which also covers case (B), to keep the paper self-contained.
\end{remark}

\subsection{Proof of part (b)}
Let $q\ge5$. By Theorem~\ref{param}, case (A) holds, so $u$ is even. Suppose, for a contradiction, that every $p_i$ divides $v$, that is, $V=N_0$. Then $W=1$, and $N/v=q^{\min(m_0,1)}\in\{1,q\}$; if $N/v=q$ then $q\mid v$.

\emph{The Lehmer pair.} Put $\alpha=u\sqrt C+\mu vi$ and $\beta=u\sqrt C-\mu vi$. Then
\[
(\alpha+\beta)^2=4Cu^2,\qquad \alpha\beta=Cu^2+v^2=y,\qquad (\alpha-\beta)^2=-4v^2 .
\]
Since $\alpha,\beta$ are the roots of $X^2-2u\sqrt C\,X+y$, whose coefficients are algebraic integers, $\alpha$ and $\beta$ are algebraic integers. Moreover $\gcd(4Cu^2,y)=1$ because $y$ is odd and $\gcd(Cu,v)=1$. Finally $\alpha/\beta=\alpha^2/y=(Cu^2-v^2+2\mu uv\sqrt{-C})/y$ lies in $K$, and it is not $\pm1$ because $uv\ne0$; since the only roots of unity in $K$ are $\pm1$, $(\alpha,\beta)$ is a Lehmer pair. Subtracting from $x\sqrt C+Ni=\delta\alpha^q$ its complex conjugate gives $\alpha^q-\beta^q=2\delta Ni$, and since $\alpha-\beta=2\mu vi$,
\begin{equation}\label{Lq}
\Lh_q(\alpha,\beta)=\delta\mu\,\frac{N}{v}\in\{\pm1,\pm q\}.
\end{equation}
If $\Lh_q(\alpha,\beta)=\pm1$, it has no prime divisor; if $\Lh_q(\alpha,\beta)=\pm q$, its only prime divisor $q$ divides $v$ and hence $(\alpha^2-\beta^2)^2=-16Cu^2v^2$. In both cases $(\alpha,\beta)$ is $q$-defective.

\emph{The case $q\ge7$.} By Theorem~\ref{bhv}, $q\in\{7,13\}$ and $|a|=4Cu^2$ belongs to $\{1,3,5,13,14\}$. This is impossible, since $u\ge2$ and $C\ge5$ give $4Cu^2\ge80$.

\emph{The case $q=5$.} With $A=(\alpha+\beta)^2=4Cu^2$ and $B=\alpha\beta=y$ we have $\Lh_5(\alpha,\beta)=\alpha^4+\alpha^3\beta+\alpha^2\beta^2+\alpha\beta^3+\beta^4=A^2-3AB+B^2$. With $w=Cu^2$ and $B=w+v^2$ this becomes
\begin{equation}\label{L5}
\Lh_5(\alpha,\beta)=5w^2-10wv^2+v^4=5X^2-4v^4,\qquad X=w-v^2,
\end{equation}
and $X$ is odd. By \eqref{Lq}, $\Lh_5(\alpha,\beta)\in\{\pm1,\pm5\}$.

If $5X^2-4v^4=\pm1$, then $(4v^2)^2-5(2X)^2=\mp4$, so $4v^2$ is a Lucas number by fact (i), and $4v^2=4$ by fact (iii). Thus $v=1$, and $N=v\cdot N/v=1$ by \eqref{Lq}, which is absurd since $N\ge3$.

If $5X^2-4v^4=\pm5$, then $N/v=5$, so $5\mid v$. Put $Z=2v^2/5$, a positive integer. Dividing $4v^4=5X^2\mp5$ by $5$ gives $5Z^2=X^2\mp1$, that is, $X^2-5Z^2=\pm1$. Since $X$ is odd, $X\ne0$, and fact (ii) applied to $(|X|,Z)$ gives $Z=F_{3n}/2$ for some $n\ge1$, so
\[
F_{3n}=2Z=\frac{4v^2}{5}=5\Bigl(\frac{2v}{5}\Bigr)^2 .
\]
By fact (iv), $3n=5$, which is impossible.

Hence $V\ne N_0$, that is, some $p_i$ does not divide $v$. \qed

\subsection{Proof of Corollary~\ref{sieve}}
By Theorem~\ref{primdiv}(b), some $p_i$ does not divide $v$, and by part (a) it satisfies $p_i\equiv\leg{-C}{p_i}\pmod q$; this gives the first assertion. If $q\ge p_i+2$ for all $i$, then $q\nmid p_i\pm1$, because $0<p_i-1<p_i+1<q$; so no $p_i$ satisfies the congruence. Finally, a prime $p_i\ne q$ with $p_i\not\equiv\leg{-C}{p_i}\pmod q$ must divide $v$ by part (a), and then $p_i^{m_i}\mid v$ by Definition~\ref{admissible}. \qed

%======================================================================
\section{The exponent \texorpdfstring{$q=3$}{q=3}}\label{sec:three}
%======================================================================

\subsection{Proof of Theorem~\ref{qthree}}
Let $u,v$ be as in Theorem~\ref{param} with $q=3$. Since $v$ is $3$-admissible, $m_0\ne1$, and $v=V$, $N/v=W$ if $m_0=0$, while $v=3^{m_0-1}V$, $N/v=3W$ if $m_0\ge2$, where $VW=N_0$ and $\gcd(V,W)=1$. Moreover $\gcd(Cu,V)=1$ since $\gcd(Cu,v)=1$, and $V,W$ are odd. Recall $\R(C,u,v,3)=Cu^2-3v^2$ and $\I(C,u,v,3)=3Cu^2-v^2$.

\emph{Case (A), $m_0=0$.} Here $W=|3Cu^2-V^2|$, so $3Cu^2=V^2+\varepsilon W$ for a sign $\varepsilon$. As $u$ is even, $12\mid3Cu^2$; as $V$ is odd and $3\nmid V$, $V^2\equiv1\pmod{12}$. Hence $\varepsilon W\equiv-1\pmod{12}$. The formulas for $x,y$ are those of (A). This is (A0).

\emph{Case (A), $m_0\ge2$.} Here $3W=|3Cu^2-3^{2m_0-2}V^2|$, so $Cu^2=3^{2m_0-3}V^2+\varepsilon W$. Since $4\mid Cu^2$, $3^{2m_0-3}\equiv3\pmod4$ and $V^2\equiv1\pmod4$, we get $\varepsilon W\equiv1\pmod4$. With $v^2=3^{2m_0-2}V^2$, (A) gives $x=u|Cu^2-3^{2m_0-1}V^2|$ and $y=Cu^2+3^{2m_0-2}V^2$. This is (A3).

\emph{Case (B), $m_0=0$.} Here $8W=|3Cu^2-V^2|$, so $3Cu^2=V^2+8\varepsilon W$. Reducing modulo $3$ and using $V^2\equiv1$ gives $8\varepsilon W\equiv-1$, that is, $\varepsilon W\equiv1\pmod3$. This is (B0).

\emph{Case (B), $m_0\ge2$.} Here $24W=|3Cu^2-3^{2m_0-2}V^2|$, so $Cu^2=3^{2m_0-3}V^2+8\varepsilon W$, and (B) gives the stated $x,y$. This is (B3).

Exactly one case occurs by Theorem~\ref{param}. The prime factors of $W$ are exactly the $p_i$ dividing $N/v$ but not $v$, so they satisfy $\ell\equiv\leg{-C}{\ell}\pmod3$ by Theorem~\ref{primdiv}(a).

For the converse, put $v=V$ if $m_0=0$ and $v=3^{m_0-1}V$ if $m_0\ge2$. We claim that $\gcd(Cu,v)=1$. This is the hypothesis $\gcd(Cu,V)=1$ when $m_0=0$. When $m_0\ge2$, reducing the displayed equation of (A3) or (B3) modulo $3$ gives $Cu^2\equiv\varepsilon W$ or $Cu^2\equiv 8\varepsilon W\pmod 3$, so $3\nmid Cu$ because $3\nmid W$; together with $\gcd(Cu,V)=1$ this proves the claim. Next, the displayed equation says that $N:=3^{m_0}VW$ equals $v\,|\I(C,u,v,3)|$ in (A0) and (A3), and $v\,|3Cu^2-v^2|/8$ in (B0) and (B3); and the displayed $x,y$ are those of (A) or (B). Since $N$ is odd, the converse of Theorem~\ref{param} applies and completes the proof. \qed

\subsection{Proof of Corollary~\ref{pure}}
Suppose first that $v=1$. In case (A), $N=|\I(C,u,1,q)|$, and the formulas of (A) give (I). In case (B), $8N=|3Cu^2-1|=3Cu^2-1$, and (B) gives $8x=u(Cu^2-3)$ (note $Cu^2\ge11$) and $4y=Cu^2+1$; this is (I$'$). Conversely (I) and (I$'$) are (A) and (B) with $v=1$.

Suppose now that $W=1$; then $q=3$ by Theorem~\ref{primdiv}(b). We put $A=Cu^2$.

If $m_0=0$ we have $V=N$. In case (A), $\varepsilon\equiv-1\pmod{12}$, so $\varepsilon=-1$ and $N^2=3A+1$; then $y=A+N^2=4A+1$ and $\R(C,u,N,3)=A-3(3A+1)=-(8A+3)$, so $x=u(8A+3)$. This is (II). In case (B), $\varepsilon\equiv1\pmod3$, so $\varepsilon=1$ and $N^2=3A-8$; then $4y=A+N^2=4A-8$ and $8x=u|A-3(3A-8)|=8u(A-3)$. This is (II$'$).

If $m_0\ge2$ we have $v=N/3$ and $3^{2m_0-3}V^2=N^2/27$. In case (A), $\varepsilon\equiv1\pmod4$, so $\varepsilon=1$ and $N^2=27(A-1)$; then $v^2=N^2/9=3(A-1)$, $y=A+v^2=4A-3$ and $\R(C,u,v,3)=A-9(A-1)=-(8A-9)$, so $x=u(8A-9)$, where $8A-9>0$ as $A\ge20$. This is (III). In case (B), $N^2=27(A-8\varepsilon)$, $v^2=3(A-8\varepsilon)$, $4y=A+v^2=4A-24\varepsilon$ and $8x=u|A-9(A-8\varepsilon)|=8u(A-9\varepsilon)$, where $A-9\varepsilon>0$ as $A\ge11$. This is (III$'$).

Conversely, each of (II)--(III$'$) is (A) or (B) with $v=N$ or $v=N/3$, and the converse of Theorem~\ref{param} applies once we check $\gcd(Cu,v)=1$. In (II) and (II$'$), a common prime factor of $Cu$ and $N$ divides $1$, respectively $8$, which is impossible as $N$ is odd. In (III) and (III$'$), $v^2=3(Cu^2-\eta)$ with $\eta=1$, respectively $\eta=8\varepsilon$; a common prime factor $\ell$ of $Cu$ and $v$ divides $3\eta$, so $\ell=3$; but then $3\nmid Cu^2-\eta$ and the $3$-adic valuation of $v^2$ would be $1$, which is absurd.

Finally let $k=1$, $N=p_1^{m_1}$. If $p_1\ne q$, then $N_0=N$ and $V\in\{1,N\}$, so $v=1$ or $W=1$. If $p_1=q$, then $N_0=1$ and $W=1$. Hence every solution satisfies one of (I)--(III$'$), and exactly one by the uniqueness in Theorem~\ref{param}. If $q\ge5$, then $W\ne1$ by Theorem~\ref{primdiv}(b), so $v=1$, family (I) occurs, and $q\nmid N$ by Lemma~\ref{rigid}(iii). \qed

\begin{remark}\label{SSremark}
For $k=1$, Corollary~\ref{pure} is the classification of \cite{SS}: the conditions there are recovered with $N=p^m$; in (III) and (III$'$) one has $N=3^m$ with $m\ge2$, and $N^2=27(Cu^2-1)$ is $Cu^2=3^{2m-3}+1$. For $k\ge2$ the statements (I)--(III$'$) still describe all solutions with $v=1$ or $W=1$, and the remaining (mixed) solutions are given by Theorem~\ref{qthree} for $q=3$ and by Theorem~\ref{param}(A) for $q\ge5$.
\end{remark}

%======================================================================
\section{Examples and a computational check}\label{sec:examples}
%======================================================================

The following table lists examples with at least two primes $p_i$ for each type. In each row the hypotheses of \eqref{maineq} hold: $C>3$ is odd and square-free, $p_i\nmid C$, $q\nmid h(-C)$, $\gcd(x,y)=1$ and $y$ is odd. Class numbers were computed by counting reduced binary quadratic forms of discriminant $-C$ or $-4C$.

\begin{center}
\renewcommand{\arraystretch}{1.25}
\begin{tabular}{@{}lllll@{}}
\toprule
Type & identity & $u$ & $v$ & $h(-C)$\\
\midrule
(I), $q=3$ & $13\cdot98^2+(5\cdot31)^2=53^3$ & $2$ & $1$ & $2$\\
(I), $q=5$ & $7\cdot1018^2+(11\cdot331)^2=29^5$ & $2$ & $1$ & $1$\\
(I), $q=7$ & $5\cdot586^2+(13^2\cdot251)^2=21^7$ & $2$ & $1$ & $2$\\
(II) & $7\cdot96804^2+(5\cdot11)^2=4033^3$ & $12$ & $55$ & $1$\\
(III) & $19\cdot1198^2+(3^2\cdot5)^2=301^3$ & $2$ & $15$ & $1$\\
(I$'$) & $227\cdot28^2+(5\cdot17)^2=57^3$ & $1$ & $1$ & $5$\\
(II$'$) & $1411\cdot1408^2+(5\cdot13)^2=1409^3$ & $1$ & $65$ & $4$\\
(III$'$), $\varepsilon=-1$ & $67\cdot76^2+(3^2\cdot5)^2=73^3$ & $1$ & $15$ & $1$\\
(III$'$), $\varepsilon=1$ & $155\cdot146^2+(3^2\cdot7)^2=149^3$ & $1$ & $21$ & $4$\\
\midrule
(A0) & $7\cdot94^2+(5\cdot59)^2=53^3$ & $2$ & $5$ & $1$\\
(A3) & $5\cdot14^2+(3^2\cdot17)^2=29^3$ & $2$ & $3$ & $2$\\
(B0) & $43\cdot4^2+(5\cdot13)^2=17^3$ & $1$ & $5$ & $1$\\
(B3) & $43\cdot2^2+(3^2\cdot5)^2=13^3$ & $1$ & $3$ & $1$\\
(A), $q=5$ & $5\cdot1990^2+(3\cdot281)^2=29^5$ & $2$ & $3$ & $2$\\
(A), $q=5$ & $13\cdot14342^2+(5^2\cdot229)^2=77^5$ & $2$ & $5$ & $2$\\
(A), $q=7$ & $5\cdot32006^2+(3\cdot36709)^2=29^7$ & $2$ & $3$ & $2$\\
\bottomrule
\end{tabular}
\end{center}

The rows below the middle rule are mixed solutions ($v>1$ and $N/v$ has a prime factor not dividing $v$). In the row $13\cdot14342^2+(5^2\cdot229)^2=77^5$ the prime $q=5$ divides $N$: as predicted by Theorem~\ref{param}, $25\mid N$, $v=5$ and $N/v=5\cdot229$. Theorem~\ref{primdiv}(a) can be checked on each row; for instance $\leg{-7}{59}=-1$ and $3\mid60$, $\leg{-5}{281}=1$ and $5\mid280$, $\leg{-13}{229}=-1$ and $5\mid230$.

\begin{remark}\label{computation}
As a check on Theorems~\ref{param} and~\ref{primdiv}, we enumerated all solutions of $Cx^2+N^2=y^q$ with $N>1$ odd, $\gcd(x,y)=1$ and $y$ odd, for odd square-free $5\le C\le199$ and $q=3,5,7,11,13$ with $y\le600,90,25,8,6$, respectively, by testing in exact integer arithmetic whether $y^q-Cx^2$ is a perfect square. Of the $1271$ solutions found, $752$ satisfy $\gcd(N,C)=1$ and $q\nmid h(-C)$ ($731$ with $q=3$, $20$ with $q=5$, one with $q=7$), and $686$ of these have at least two distinct prime factors of $N$. For each of the $752$ solutions we found exactly one pair $(u,v)$ as in (A) or (B), checked that $v$ is $q$-admissible, that every prime of $N$ not dividing $v$ satisfies the congruence of Theorem~\ref{primdiv}(a), and that for $q\ge5$ such a prime exists. No exception occurred. Separately, we checked the converse of Theorem~\ref{param} for all odd square-free $5\le C\le299$: in case (A) for $q\in\{3,5,7,11\}$, even $u\le12$ and odd $v<60$, and in case (B) for odd $u\le13$ and odd $v<60$, always with $\gcd(Cu,v)=1$ ($68597$ cases in all).
\end{remark}

\begin{remark}\label{yeven}
The hypothesis that $y$ is odd cannot be omitted. By Lemma~\ref{yodd} it matters only when $C\equiv7\pmod 8$, and then there are solutions with $y$ even, for example
\[
71\cdot1^2+(3\cdot7)^2=8^3,\qquad 7\cdot11^2+(3\cdot19)^2=16^3,
\]
with $h(-71)=7$ and $h(-7)=1$. They violate the conclusions of Theorem~\ref{param}: $x$ is odd and $3$ divides $N$ exactly once.
\end{remark}

%======================================================================
\section{Concluding remarks}\label{sec:final}
%======================================================================

\textbf{Composite exponents.} If $n$ is odd and $q\mid n$ is a prime with $q\nmid h(-C)$, a solution of $Cx^2+N^2=y^n$ gives the solution $(x,y^{n/q})$ of \eqref{maineq}. Hence Theorem~\ref{param} applies with $y$ replaced by $y^{n/q}$; in particular Corollary~\ref{sieve} shows that $Cx^2+p_1^{2m_1}\cdots p_k^{2m_k}=y^n$ has no solution with $y$ odd when $\gcd(n,h(-C))=1$ and $n$ has a prime factor $q\ge\max(5,p_1+2,\dots,p_k+2)$.

\textbf{Open questions.} For fixed $C$ and fixed primes $p_i$, each admissible distribution of the primes between $v$ and $N/v$ leads to an equation of the form $v|\I(C,u,v,q)|=N$ in the unknowns $u$ and $m_i$; for $q=3$ these are the auxiliary equations of Theorem~\ref{qthree}. It would be interesting to solve these completely for small sets of primes, for instance by the methods of \cite{Tho}, to describe the solutions with $y$ even when $C\equiv7\pmod8$, and to treat the case $q\mid h(-C)$.

\subsection*{Funding}
The authors did not receive support from any organization for the submitted work.

\subsection*{Data availability}
No data sets were generated or analysed in this study, apart from the computations described in Remark~\ref{computation}; the programs used are available from the authors on request.

\subsection*{Conflict of interest}
The authors declare that they have no conflict of interest.

\end{document}